\documentclass[12pt]{article} 

\usepackage{graphicx} 
\usepackage{amsmath}
\usepackage{amsfonts}
\usepackage{apacite}

\newcommand {\omegavec}{{\boldsymbol{\omega}}}

\newcommand {\lambdavec}{{\boldsymbol{\lambda}}}

\newcommand {\zetavec}{\boldsymbol{\zeta}}

\newcommand {\thetavec}{{\boldsymbol{\theta}}}

\newfont{\pseudocode}{cmtt10}

\newtheorem{thm}{Theorem}

\newtheorem{lem}{Lemma}

\newtheorem{assume}{Assumption}

\newtheorem{rem}{Remark}
\newtheorem{problem}{Problem}

\newcommand{\bfx}{\mathbf{x}}
\newcommand{\bff}{\mathbf{f}}
\newcommand{\bfb}{\mathbf{b}}
\newcommand{\bfc}{\mathbf{c}}

\newcommand{\bfw}{\mathbf{w}}

\newcommand{\bfu}{\mathbf{u}}
\newcommand{\bfg}{\mathbf{g}}

\newcommand{\bfq}{\mathbf{q}}

\newcommand{\pd}{{\partial}}

\newcommand{\Real}{\mathbb{R}}
\newcommand{\Natural}{\mathbb{N}}
\newcommand{\Numbers}{\mathbb{Z}}

\newcommand{\norm}[1]{\left\Vert#1\right\Vert}

\newcommand {\ess}{{\mathrm{ess}}}
\newcommand{\norms}[1]{\norm{#1}_{\mathcal{A}}}
\newcommand{\normstwo}[2]{\norm{#1}_{#2}}

\begin{document}

\baselineskip 0.7cm

\title{\bf {Adaptive classification of temporal signals in fixed-weights recurrent neural networks: an existence proof}}

\author{Ivan Y. Tyukin\thanks{Department of Mathematics, School of Mathematics and Computer Science, University of Leicester,  University Rd., Leicester, LE1 7RH,
        UK (e-mail: I.Tyukin@le.ac.uk); and Laboratory for Perceptual Dynamics, Brain Science Institute, RIKEN (Institute for Physical and Chemical
 Research), Wako-shi, Saitama, 351-0198, Japan (e-mail: tyukinivan@brain.riken.jp)}, Danil Prokhorov\thanks{Toyota Technical Center, Ann Arbor 48105, USA (e-mail: dvprokhorov@gmail.com)}, Cees van Leeuwen
 \thanks{Laboratory for Perceptual Dynamics, Brain Science Institute, RIKEN (Institute for Physical and Chemical
 Research), Wako-shi, Saitama, 351-0198, Japan (e-mail: ceesvl@brain.riken.jp)}
}

\date{}

\maketitle

\begin{abstract}

We address  the important theoretical question why a recurrent
neural network with fixed weights can adaptively classify
time-varied signals in the presence of additive noise and parametric
perturbations.  We provide a mathematical proof assuming that
unknown parameters are allowed to enter the signal nonlinearly and
the noise amplitude is sufficiently small.
\end{abstract}


\section{Introduction}
\label{Introduction}

Recurrent Neural Networks (RNN) with fixed weights are known to be
able to solve
problems of adaptive classification, recognition, and control
\shortcite{Conf:IJCNN-02,Conf:ICNN-96:Feldkamp,Conf:ICNN-97:Feldkamp,IEEE_TNN:1999:Younger,Conf:IJCNN-2001:Lo}.
When the objects to be classified are static, e.g. still images or
vectors in $\Real^n$, the way the fixed-weight RNN solves problems
is usually characterized in terms of convergence of the RNN state to
an attractor \shortcite{Hopfield82,Fuchs88}. Each attractor
corresponds to a specific class of objects and its basin determines
which objects belong to the class. Conditions specifying convergence
to an attractor are widely available in this case,
\shortcite{IEEE_TSMC:Cohen:1983,IEEE_TCS:Michel:1990,IEEE_TNN:Yang:1994,IEEE_TNN:Amari:2001,NeuralComp:Chen:2003}
to name a few.
%

When the objects to be classified are dynamic, for instance
nonlinearly parameterized functions of time of which the  parameters
are unknown a-priori, no adequate theory exists that explains why
the fixed-weight RNN approach is successful. At present, theoretical
results are available to demonstrate that a single fixed-weight RNN
of a certain type can {\it approximate} the solutions of multiple
dynamical systems \shortcite{Neural_Computation:Chen:2002}. Hence in
principle, a fixed-weight RNN can behave adaptively with respect to
changes of its input signals. These theoretical results, however,
are restricted to the class of parameter replacement networks
\shortcite{IEEE_TNN:Chen:1995}. The structure of these networks
differs from that of the more commonly used recurrent multilayered
perceptrons. Whether adaptive behavior is inherent to other types of
RNN, therefore, remains an unresolved theoretical issue. In spite of
plausibility arguments given by several authors
\shortcite{Conf:ICNN-97:Feldkamp,Conf:IJCNN-02}, no formal proof has
been made
available, to the best of our knowledge. 

%

In this paper we consider adaptive behavior in fixed-weight RNNs
from the standpoint of their ability to {\it classify temporal
signals} adaptively. We provide a formal proof that {\it
continuous-time} recurrent neural networks with fixed weighs can
successfully classify and recognize nonlinear functions of time and
unknown parameter. These functions are allowed to be nonlinearly
parameterized.  The main idea behind our results consists of
presenting a prototype dynamical system which solves the recognition
problem. This is followed by a proof that a RNN with fixed weights
can realize this system. We construct such a system using the
concepts of relaxation times and weakly attracting sets
\shortcite{CommMathPhys:Milnor:1985,Gorban:2004} as well as the
tests for convergence to such sets obtained in our earlier work
\shortcite{ArXive:Non-uniform:2006}. To show that our system can
indeed be realized by a RNN with fixed weights we employ classical
results on function approximation by feed-forward networks
\shortcite{Cybenko}.

The paper is organized as follows. Section \ref{Notation} describes
notational agreements. In Section \ref{Problem Formulation} we
provide a mathematical statement of the problem, Section
\ref{Main_Results} contains the main results, and Section
\ref{Conclusion} concludes the paper.

\section{Notational Preliminaries}\label{Notation}

\hspace{5mm} $\bullet$ Symbol $\Real$ defines the field of real
numbers, and symbol $\Real_{\geq c}$, $c\in\Real$ stands for the
following set $\Real_{\geq c}=\{x\in\Real|x\geq c\}$, and $\Real_{>
c}=\{x\in\Real|x> c\}$.

$\bullet$ Symbol $\Real^n$ stands for an $n$-dimensional linear
space over the field of reals.

$\bullet$ $\mathcal{C}^k$ denotes the space of functions that are at
least $k$ times differentiable.

$\bullet$ Symbol $\mathcal{K}$ denotes the class of all strictly
increasing functions $\kappa: \Real_{\geq 0}\rightarrow \Real_{\geq
0}$ such that $\kappa(0)=0$; symbol $\mathcal{K}_\infty$ denotes the
class of all functions $\kappa\in\mathcal{K}$ such that
$\lim_{s\rightarrow\infty}\kappa(s)=\infty$.

$\bullet$ Symbol $\oplus$ denotes concatenation of two vectors.

$\bullet$ The solution of a system of differential equations
$\dot{\bfx}=\bff(t,\bfx,\thetavec,\bfu(t))$,
$\bff:\Real\times\Real^n\times\Real^d\times\Real^m\rightarrow\Real^n$,
$\bff\in\mathcal{C}^0$, $\bfu:\Real_{\geq 0}\rightarrow\Real^m$,
$\thetavec\in\Real^d$ passing through point $\bfx_0$ at $t=t_0$ will
be denoted for $t\geq t_0$ as $\bfx(t,\bfx_0,t_0,\thetavec,\bfu)$,
or simply as $\bfx(t,\bfx_0)$ or $\bfx(t)$ if it is clear from the
context what the values of $\bfx_0,\thetavec$ are and how the
function $\bfu(t)$ is defined.

$\bullet$ By ${L}^n_\infty[t_0,T]$, $t_0\geq 0$, $T\geq t_0$ we
denote the space of all functions $\bff:\Real_{\geq
0}\rightarrow\Real^n$ such that $\|\bff\|_{\infty,[t_0,T]}=\ess
\sup\{\|\bff(t)\|,t \in [t_0,T]\}<\infty$;
$\|\bff\|_{\infty,[t_0,T]}$ stands for the ${L}^n_\infty[t_0,T]$
norm of $\bff(t)$.

$\bullet$ Let $\mathcal{A}$ be a set in ${\Real^n}$ and $\|\cdot\|$
be the usual Euclidean norm in $\Real^n$. By the symbol
$\norms{\cdot}$ we denote the following induced norm:
\[
\norms{\bfx}=\inf_{\bfq\in\mathcal{A}}\{\|\bfx-\bfq\|\}
\]
In case $x$ is a scalar and $\Delta\in\Real_{>0}$, notation
$\|x\|_\Delta$ stands for the following
\[
\|x\|_{\Delta}=\left\{
                    \begin{array}{ll}
                    |x|-\Delta, & |x|> \Delta\\
                    0, & |x|\leq \Delta
                    \end{array}
                    \right.
\]

\section{Problem Formulation}\label{Problem Formulation}

Consider the following set of signals
\begin{equation}\label{eq:signals}
\begin{split}
\mathcal{F}&=\{f_i(\xi(t),\theta_i)\},  \  i\in\{1,\dots,N_f\}, \\
& f_i: \Real\times\Real\rightarrow\Real, \
f_i(\cdot,\cdot)\in\mathcal{C}^0,\\
& \xi:\Real_{\geq 0}\rightarrow\Real, \
\xi(\cdot)\in\mathcal{C}^1\cap L_{\infty}[0,\infty]
\end{split}
\end{equation}
where $\theta_i\in\Omega_\theta\subset\Real$ are parameters of which
the values are unknown a-priori,
$\Omega_\theta=[\theta_{\min},\theta_{\max}]$ is a bounded interval,
and $\xi(t)$ is a known and bounded function. Signals
$f_i(\xi(t),\theta_i)$ represent relevant physical variables of an
object.

For the given functions $f_i(\xi(t),\theta_i)$ and $\xi(t)$ we say
that $\theta_i$ is {\it equivalent} to $\theta_i'$ iff
\begin{equation}\label{eq:equivalece_notion}
f_i(\xi(t),\theta_i)=f_i(\xi(t),\theta_i') \ \forall \
t\in\Real_{\geq 0}.
\end{equation}
Hence an equivalence class for $\theta_i\in\Omega_\theta$  can be
defined as
\begin{equation}\label{eq:equivalent_class}
E_i(\theta_i)=\{\theta_i'\in\Real| \
f_i(\xi(t),\theta_i)=f_i(\xi(t),\theta_i') \ \forall \
t\in\Real_{\geq 0}\}
\end{equation}
Equivalence classes (\ref{eq:equivalent_class}) determine sets of
indistinguishable parameterizations of the $i$-th signal. It is
natural, therefore, to restrict ourselves to the problem of
recognizing signals (\ref{eq:signals}) up to their equivalence
classes.

With respect to the equivalence classes $E_i(\theta_i)$, we further
assume that there is at least one point $\theta_{0}\in\Real$ such
that
\begin{equation}\label{eq:non-degeneracy_recognition}
\|\theta_{0}\|_{E_i(\theta_i)}\geq \Delta_\theta\in\Real_{>0} \
\forall \ \theta_i\in\Omega_\theta.
\end{equation}
Requirement (\ref{eq:non-degeneracy_recognition}) is a technical
assumption. It holds, however, for a wide range of practically
relevant situations in which the union of $E_i(\theta_i)$ for all
$i$ and $\theta_i$ belongs to an interval of $\Real$. Furthermore,
it allows us to exclude from consideration pathological cases in
which almost all points in $\Omega_\theta$ are indistinguishable in
the sense of condition (\ref{eq:equivalece_notion}).

In many systems, artificial or natural, measured physical
quantities, represented here by  signals $f_i(\xi(t),\theta_i)$, are
often unavailable. This is because a measurement device is involved
in measuring $f_i(\xi(t),\theta_i)$. Given that signals
$f_i(\xi(t),\theta_i)$ are functions of time, inherent dynamical
properties of a measurement device would distort the measured
values. Our present study takes this possibility into account. To do
so we consider the case where signals $f_i(\xi(t),\theta_i)$ are
affected by additive bounded noise and pass through nonlinear
filters with uncertain dynamics. In particular, we assume that
instead of functions $f_i(\xi(t),\theta_i)$ we access variables
$s_i(t,s_{i,0},\theta_i,\eta_i(t))$, which are solutions to the
following ordinary differential equation:
\begin{equation}\label{eq:filtered_signals}
\begin{split}
\dot{s}_i&=-\varphi_i(s_i)+f_i(\xi(t),\theta_i) + \eta_i(t),  \\
&s_i(t_0)=s_{i,0}, \ s_{i,0}\in\Omega_s\subset\Real
\end{split}
\end{equation}
In (\ref{eq:filtered_signals}) the function
$\eta_i:\Real_{>0}\rightarrow\Real$:
\begin{equation}\label{eq:additive_noise}
\eta_i(t)\in L_{\infty}[0,\infty], \
\|\eta_i(t)\|_{\infty,[0,\infty]}\leq \Delta_{\eta}\in\Real_{\geq 0}
\end{equation}
corresponds to measurement noise. The value of $\Delta_{\eta}$ in
(\ref{eq:filtered_signals}) is supposed to be known, while the
values of initial conditions $s_i(t_0)$ and functions
$\varphi_i:\Real\rightarrow\Real, \ \varphi(\cdot)\in\mathcal{C}^1$
in (\ref{eq:filtered_signals}) are assumed to be uncertain. We do,
however, require that $\Omega_s=[s_{\min},s_{\max}]$ is an interval
and that  the functions $\varphi_i(s_i)$ satisfy the following
constraint:
\begin{equation}\label{eq:filter_stability}
 \forall s_i\in\Real \Rightarrow  \
\varphi_{\min}\leq \frac{\pd \varphi_i(s_i)}{\pd s_i}\leq
\varphi_{\max}, \ \varphi_{\min},\varphi_{\max}\in\Real_{>0}.
\end{equation}
Condition (\ref{eq:filter_stability}) ensures that filters
(\ref{eq:filtered_signals}) are convergent \shortcite{Pavlov:2004},
e.g. the dynamics of each  variable
$s_i(t,s_{i,0},\theta_i,\eta_i(t))$ at $t\rightarrow\infty$  is
uniquely determined in the absence of noise by
$f_i(\xi(t),\theta_i)$, and the effects of initial conditions
$s_{i,0}$ vanish with time asymptotically.

A recurrent neural network is defined by the following set of
differential equations:
\begin{eqnarray}\label{eq:RNN_recurrent}
& & \dot{x}_j= \sum_{m=1}^{N} c_{j,m}\sigma(\bfw^{T}_{j,m}(\xi(t)\oplus s(t) \oplus\bfx)+b_{j,m}), \ j\in\{1,\dots,N_x\},\\
& & \bfx=\mathrm{col}(x_1,\dots,x_{N_x}), \
\bfx(t_0)=\bfx_0,\nonumber
\end{eqnarray}
where functions $\sigma:\Real\rightarrow\Real$ are sigmoid. Vectors
$\bfc_j=\mathrm{col}(c_{j,1},\dots,c_{j,N})$,
$\bfb_j=\mathrm{col}(b_{j,1},\dots,b_{j,N})$ and matrices
$\boldsymbol{W}_j=(\bfw_{j,1},\dots,\bfw_{j,N})$ are the RNN
parameters. Functions $\xi(t)$, $s(t):\Real_{\geq
0}\rightarrow\Real$, $\xi(t), s(t)\in\mathcal{C}^0$ are inputs;
$\bfx$ is the state vector, and $\bfx_0$ is a vector of initial
conditions.
%
%

According to notation (\ref{eq:RNN_recurrent}) the network maps two
functions of time $\xi(t)$, $s(t)$ into the functions
$x_1(t,\bfx_0),\dots, x_{N_x}(t,\bfx_0)$, which are the solutions of
(\ref{eq:RNN_recurrent}). In what follows we will consider variables
$\xi(t)$, $s(t)$ as inputs to the network. While the variable
$\xi(t)$ is known a-priori, variable $s(t)$ is allowed to vary
within the set of functions $s_i(t,s_{i,0},\theta_i,\eta_i(t))$,
which are the solutions of (\ref{eq:filtered_signals}). In
particular, we assume that the following condition is satisfied:
\begin{assume}[Existence]\label{assume:signal} There exist $i\in N_f$,
$\theta_i\in\Omega_\theta$, $s_{i,0}\in\Omega_s$ and $\eta_i(t)$
specified by (\ref{eq:additive_noise}) such that
\begin{equation}\label{eq:class_definition}
s(t)=s_i(t,s_{i,0},\theta_i,\eta_i(t)) \ \forall \ t\geq 0.
\end{equation}
\end{assume}
We aim to determine if there is a network of type
(\ref{eq:RNN_recurrent}) which is able to recover uncertain
parameters $i$ and $\theta_i$ from the input
 $s(t)$\footnote{Because filters (\ref{eq:filtered_signals}) are convergent, the effect of uncertainty in parameter $s_{i,0}$ vanishes with time
 exponentially. Hence the only effective uncertainties are $i$ and $\theta_i$.}, $t\geq t_0\in\Real_{\geq0}$ within a finite  interval of
 time 
 for all $t_0\in\Real_{\geq 0}$ . Informally, this means
 that there exist two sets of functions of network state $\bfx$ and input $s(t)$:
\begin{equation}\label{eq:RNN_output}
\begin{split}
&\{h_{f,j}(\bfx(t),s(t))\}, \  \{h_{\theta,j}(\bfx(t),s(t))\}, \\
& h_{f,j}:\Real^{N_x}\times\Real\rightarrow\Real, \
h_{\theta,j}:\Real^{N_x}\times\Real\rightarrow\Real, \
j\in\{1,\dots,N_f\},
\end{split}
\end{equation}
 such that the values of $i$ and $\theta_i$ can be inferred from $\{h_{f,j}(\bfx(t),s(t))\}$,
 $\{h_{\theta,j}(\bfx(t),s(t))\}$ respectively within a given finite interval of time. Formally we can
state this as follows:

\begin{problem}\label{problem:RNN_existence} Consider class $\mathcal{F}$
of signals (\ref{eq:signals}), where the function $\xi(t)$ is known,
and the values of parameters $\theta_i$ are unknown a-priori.
Determine a recurrent neural network (\ref{eq:RNN_recurrent}) such
that the following properties hold:

1) there is a set of initial conditions $\Omega_x$ such that
$\bfx(t,\bfx_0)$ is bounded for all $\bfx_0\in\Omega_x$ and $t\geq
t_0\in\Real_{\geq 0}$; the volume of $\Omega_x$ is nonzero;

2) there exists a set of output functions (\ref{eq:RNN_output}) such
that, for all  $\theta_i\in\Omega_{\theta}$, $s_{i,0}\in\Omega_s$,
$t_0\in\Real_{\geq 0}$, $\bfx_0\in\Omega_x$, and functions
$\eta_i(t)$ given by (\ref{eq:additive_noise}), condition
(\ref{eq:class_definition})
implies existence of a constant $\mathcal{T}\in\Real_{>0}$, time
instant $t'\in(t_0,t_0+\mathcal{T})$, (arbitrarily large)
$T^\ast\in\Real_{>0}$,  and (arbitrarily small)
$\varepsilon\in\Real_{>0}$ and $\mathcal{D}\in\mathcal{K}_\infty$
such that
\begin{equation}\label{eq:class_recognition}
\begin{split}
& \| h_{f,i}(\bfx(t),s(t)) \|_{\infty,[t',t'+T^\ast]}<\varepsilon +
\mathcal{D}(\Delta_\eta),  \\
& \inf_{\theta_i'\in E(\theta_i)} \|
 h_{\theta,i}(\bfx(t),s(t))-\theta_i'\|_{\infty,[t',t'+T^\ast]}<\varepsilon +
 \mathcal{D}(\Delta_\eta).
 \end{split}\nonumber
 \end{equation}
\end{problem}

In general, this problem has no solutions for all possible
$\xi(t)\in\mathcal{C}^1$ and $f_i(\cdot,\cdot)\in\mathcal{C}^0$.
Consider, for instance, the case when
$f_i(\xi(t),\theta_i)=\sin(\xi(t)\theta_i)$ and
\[
\xi(t)=\left\{\begin{array}{ll}
                     \sin^2(\ln(t-t_0+1)), & \sin(\ln(t-t_0+1))\geq
                     0\\
                      0, & \sin(\ln(t-t_0+1))< 0
               \end{array}
\right. \ \forall \ t\geq t_0.
\]
Time intervals  when $\xi(t)=0$ are growing unboundedly with time.
Hence for any fixed $\mathcal{T}$, $T^\ast$ there will always exist
an instant $t_0'$ such that for all $t\geq t_0'$ the lengths of
intervals when $\xi(t)=0$ exceed $\mathcal{T}+T^\ast$. For all such
intervals, solutions $s_i(t,s_{i,0},
 \theta_i,\eta_i(t))$ do not depend on $\theta_i$. Hence recovery of the
 actual values of $\theta_i$ from signal $s(t)$ cannot be achieved within a fixed time interval $[t_0, t_0+\mathcal{T}+T^\ast]$ for all $t_0\geq t_0'$. In
 order to enable a solution of the classification/recognition problem above we
 must introduce an additional constraint on the functions
 $f_i(\xi(t),\theta_i)$. This should ensure that
 variation in parameter $\theta_i$ can be detected from
 the values of $f_i(\xi(t),\theta_i)$ within a finite time interval.
We therefore require that the  following property holds:

\begin{assume}[Non-degeneracy]\label{assume:NLPE} For the set of functions
$f_i(\xi(t),\theta_i)$ specified by (\ref{eq:signals}) and all
$t\geq t_0$, $\theta_i, \ \theta_i'$ there exist a constant
$T\in\Real_{>0}$ and a strictly increasing function
$\rho:\Real_{\geq 0}\rightarrow \Real_{\geq 0}$,
$\rho\in\mathcal{K}_\infty$ such that the following condition holds:
\begin{equation}\label{eq:nonlinear_persistency}
\forall \ t\geq t_0 \ \exists \ t'\in [t,t+T]: \ \
|f_i(\xi(t'),\theta_i)-f_i(\xi(t'),\theta_i')|\geq
\rho\left(\normstwo{\theta_i}{E_i(\theta_i')}\right).
\end{equation}
\end{assume}
In case the equivalence classes $E_i(\theta_i')$ consist of single
elements, e.g. when there is a unique value of $\theta_i'=\theta_i$
satisfying (\ref{eq:equivalece_notion}), condition
(\ref{eq:nonlinear_persistency}) will have a more transparent form:
\begin{equation}\label{eq:nonlinear_persistency_unique}
\forall \ t\geq t_0 \ \exists \ t'\in [t,t+T]: \ \
|f_i(\xi(t'),\theta_i)-f_i(\xi(t'),\theta_i')|\geq
\rho(|\theta_i-\theta_i'|).
\end{equation}
These conditions simply state that within a fixed time interval the
values of $\normstwo{\theta_i}{E_i(\theta_i')}$ or
$|\theta_i-\theta_i'|$ can be inferred from the differences
$f_i(\xi(t),\theta_i)-f_i(\xi(t),\theta_i')$ for all
$t\in\Real_{\geq 0}$.


In the next section we show that the solution to Problem
\ref{problem:RNN_existence} can be obtained for the class
$\mathcal{F}$ of functions $f_i(\xi(t),\theta_i)$ that are Lipschitz
in $\theta_i$. We present these results in the form of sufficient
conditions formulated in Theorem \ref{theorem:RNN_existence}.

\section{Main Results}\label{Main_Results}

As was suggested in our previous work \shortcite{Conf:IJCNN-02}, as
well as in \shortcite{IEEE_TNN:1999:Younger} the reason why RNNs
with fixed parameters (weights) demonstrate adaptive behavior could
be found in their dynamics; supposedly, it is already sufficiently
rich to have an adequate adaptation mechanism embedded into it.
Finding  a system which satisfies requirements 1), 2) in Problem
\ref{problem:RNN_existence} and which is, at the same time,
realizable by a RNN, therefore, automatically constitutes an
existence proof. This intuition, we will show, is correct. The
result is provided in Theorem \ref{theorem:RNN_existence} below.


\begin{thm}[Existence]\label{theorem:RNN_existence} Let
functions $\xi(t)$, $f_i(\xi(t),\theta_i)$ be given and defined as
in (\ref{eq:signals}), and Assumptions \ref{assume:signal},
\ref{assume:NLPE} hold. Furthermore, suppose that
$f_i(\xi(t),\theta_i)$ are (locally) Lipschitz\footnote{Property
(\ref{eq:Lipschitz}) can be understood as a generalized Lipschitz
condition. When equivalence sets $E_i(\theta_i')$ consist of single
elements the property transforms into:
$|f_i(\xi(t),\theta_i)-f_i(\xi(t),\theta_i')|\leq D_{\theta}
|\theta_i-\theta_i'|$.}:
\begin{eqnarray}
&\exists \ D_{\theta}\in\Real_{>0}:  \
|f_i(\xi(t),\theta_i)-f_i(\xi(t),\theta_i')|\leq D_{\theta}
\normstwo{\theta_i}{E_i(\theta_i')}& \ \ \forall \ t>0,  \
\theta_i,\theta_i' \label{eq:Lipschitz} \\
&\exists \ D_{\xi}\in\Real_{>0}:  \
|f_i(\xi,\theta_i)-f_i(\xi',\theta_i)|\leq D_{\xi} |\xi - \xi'| & \
\ \ \ \ \forall \ \theta_i,  \xi,  \xi' \label{eq:Lipschitz_xi}
\end{eqnarray}
and the time-derivative of $\xi(t)$ is bounded:
\begin{equation}\label{eq:xi_derivative}
 \left|\frac{d}{dt} \xi(t)\right|\leq \pd \xi_\infty \ \forall  \
 t\geq 0.
\end{equation}

Then for any $T^\ast\in\Real_{>0}$, $\varepsilon\in\Real_{>0}$ there
is a recurrent neural network (\ref{eq:RNN_recurrent}) satisfying
the requirements of Problem \ref{problem:RNN_existence}, provided
that the upper bound $\Delta_\eta$ for the
$L_\infty[0,\infty]$-norms of the disturbance terms, $\eta_i(t)$, is
sufficiently small.
\end{thm}

{\it Proof of Theorem \ref{theorem:RNN_existence}.} We prove the
theorem in four steps. First, we present a dynamical system which
will be referred to as the {\it convergence prototype}. We select
this system in the following class of differential-algebraic
equations:
\begin{eqnarray}
& &\begin{aligned}
\dot{\hat{s}}_i&=-\varphi_i(\hat{s}_i)+f_i(\xi(t),\hat{\theta}_i)
\end{aligned}\label{eq:convergence_prototype}\\
& & \begin{aligned} \hat{\theta}_i&= a+\frac{b-a}{2}(x_i+1)
\end{aligned}\label{eq:convergence_prototype:1}\\
& & \begin{aligned}
\dot{x}_i&=\gamma \|\hat{s}_i-s\|_{\varepsilon}\left(x_i-y_i-x_i(x_i^2+y_i^2)\right)\\
\dot{y}_i&=\gamma
\|\hat{s}_i-s\|_{\varepsilon}\left(x_i+y_i-y_i(x_i^2+y_i^2)\right),
\end{aligned}\label{eq:convergence_prototype:2}
\end{eqnarray}
where
\begin{equation}\label{eq:theorem:prototype_parameters}
\gamma\in\Real_{>0}, \ a,b\in \Real, \  a<\theta_{\min}, \
b>\theta_{\max}, \ \theta_0\in [a,b], \  i=1,\dots,N_f, \
\varepsilon\in\Real_{>0}.
\end{equation}
System
(\ref{eq:convergence_prototype})--(\ref{eq:convergence_prototype:2})
has a locally Lipschitz right-hand side and its solutions are
bounded for all initial conditions $\hat{s}_i(t_0)$, $x_i(t_0)$,
$y_i(t_0)\in\Real$. We show that there exist (domains  of)
$\gamma>0$, $\varepsilon>0$ and a point $\hat{s}_i(t_0)=s_0'$,
$x_i(t_0)=x_0'$, $y_i(t_0)=y_0'$, such that the trajectories passing
through this point converge to the following target set
\begin{equation}\label{eq:target_set_RNN}
\|\hat{s}_i-s_i\|_{\varepsilon}=0, \ \
\normstwo{\hat{\theta}_i}{E_i(\theta_i)}\leq
\varepsilon_{\theta}(\varepsilon).
\end{equation}
Second, we prove that there is a point $x_i(t_0)=x_0'$,
$y_i(t_0)=y_0'$ such that convergence is locally uniform with
respect to the values of uncertain $\theta_i$ and $s_{i,0}$. In
other words, for all $t_0\geq 0$, $s_{i,0}\in\Omega_s$, and
$\theta_i\in\Omega_{\theta}$ there exists $\tau>0$ such that
solutions of
(\ref{eq:convergence_prototype})--(\ref{eq:convergence_prototype:2})
with initial conditions $x_i(t_0)=x_0'$, $y_i(t_0)=y_0'$  will be in
an arbitrarily small neighborhood of (\ref{eq:target_set_RNN}) for
all $t\geq t_0+\tau$.

System
(\ref{eq:convergence_prototype})--(\ref{eq:convergence_prototype:2}),
however, is not structurally stable. That is, small perturbations of
its right-hand side might change asymptotic properties of the system
drastically. Hence, due to the inevitable approximation errors, the
chances that an RNN realization of
(\ref{eq:convergence_prototype})--(\ref{eq:convergence_prototype:2})
would solve Problem \ref{problem:RNN_existence} are slim. To
continue our argument we need to modify
(\ref{eq:convergence_prototype})--(\ref{eq:convergence_prototype:2})
such that the resulting system becomes structurally stable.

For this reason we, third, consider the perturbed version of system
(\ref{eq:convergence_prototype})--(\ref{eq:convergence_prototype:2})
\begin{eqnarray}
& & \begin{aligned}
\dot{\hat{s}}_i&=-\varphi_i(\hat{s}_i)+f_i(\xi(t),\hat{\theta}_i)\\
\hat{\theta}_i&= a+\frac{b-a}{2}(x_i+1)
\end{aligned}\label{eq:convergence_prototype_struct_stable}\\
& & \begin{aligned}
\dot{x}_i&=\gamma (\|\hat{s}_i-s\|_{\varepsilon}+\delta)\left(x_i-y_i-x_i(x_i^2+y_i^2)\right)\\
\dot{y}_i&=\gamma
(\|\hat{s}_i-s\|_{\varepsilon}+\delta)\left(x_i+y_i-y_i(x_i^2+y_i^2)\right),
\ \delta\in\Real_{>0}
\end{aligned}\label{eq:convergence_prototype_struct_stable:2}
\end{eqnarray}
aiming at achieving structural stability of an otherwise
structurally unstable system. We show that trajectories of system
(\ref{eq:convergence_prototype_struct_stable}),
(\ref{eq:convergence_prototype_struct_stable:2}) periodically visit
a small vicinity of (\ref{eq:target_set_RNN}) and stay there for an
arbitrary long time, depending on the value of $\delta$. Fourth,
given that system (\ref{eq:convergence_prototype_struct_stable}),
(\ref{eq:convergence_prototype_struct_stable:2}) is structurally
stable, we apply the results from \shortcite{Cybenko} to demonstrate
that solutions of (\ref{eq:convergence_prototype_struct_stable}),
(\ref{eq:convergence_prototype_struct_stable:2}) can be approximated
in forward time over the semi-infinite interval $[0,\infty]$ by the
state of a recurrent neural network specified by equations
(\ref{eq:RNN_recurrent}).

{\it 1. Convergence prototype.}  According to Assumption
\ref{assume:signal} there exist  $i\in\{1,\dots,N_f\}$, $s_{i,0}$,
$\theta_i$ such that $s(t)=s_i(t,s_{i,0},\theta_i,\eta_i(t))$ for
all $t\geq 0$. Consider the $i$-th subsystem of
(\ref{eq:convergence_prototype})--(\ref{eq:convergence_prototype:2})
and analyze the dynamics of the following difference:
$s_i(t)-\hat{s}_i(t)$. Denoting
\begin{equation}\label{eq:theorem:1}
\begin{split}
e_i(t)&=s(t)-\hat{s}_i(t)=s_i(t)-\hat{s}_i(t), \\
\alpha_i(t)&=\int_{0}^1\frac{\pd \varphi(s_i(t) r +
(1-r)\hat{s}_i(t))}{\pd s_i(t) r + (1-r)\hat{s}_i(t))} d r\\
\Delta
f_i(t)&=f_i(\xi(t),\theta_i)-f_i(\xi(t),\hat{\theta}_i(x_i(t)))
\end{split}
\end{equation}
and using Hadamard's lemma  we can derive the following estimate:
\begin{equation}\label{eq:theorem:2}
|e_i(t)|\leq e^{-\int_{0}^t \alpha_i(\tau)d\tau} |e_i(0)| +
\frac{1}{\varphi_{\min}} \left(1 - e^{-\int_{0}^t
\alpha_i(\tau)d\tau}\right)(\|\Delta
f_i(\tau)\|_{\infty,[0,t]}+\|\eta_i(\tau)\|_{\infty,[0,\infty]})
\end{equation}
Given that $\|\eta_i(\tau)\|_{\infty,[0,\infty]}\leq \Delta_\eta$
for all $t\geq 0$, inequality (\ref{eq:theorem:2}) implies that
\[
\left(|e_i(t)|-
\frac{\Delta_{\eta}}{\varphi_{\min}}\right)\leq e^{-\varphi_{\min}
t} \left(|e_i(0)|- \frac{\Delta_{\eta}}{\varphi_{\min}}\right) +
\frac{1}{\varphi_{\min}} \|\Delta f_i(\tau)\|_{\infty,[0,t]}
\]
Hence the following estimate holds along the trajectories of
(\ref{eq:convergence_prototype}):
\begin{equation}\label{eq:theorem:3}
\|e_i(t)\|_{\varepsilon}\leq e^{-\varphi_{\min} t}
\|e_i(0)\|_{\varepsilon} + \frac{1}{\varphi_{\min}} \|\Delta
f_i(\tau)\|_{\infty,[0,t]}, \ \
\varepsilon=\frac{\Delta_{\eta}}{\varphi_{\min}}
\end{equation}
Taking (\ref{eq:Lipschitz}), (\ref{eq:theorem:3}) into account plus
the fact that
$\normstwo{\hat{\theta}_i}{E_i(\theta_i)}=\inf_{\bar{\theta}_i\in
E_i(\theta_i)}|\hat{\theta}_i-\bar{\theta}_i|$ we can conclude that
the following inequality holds:
\begin{equation}\label{eq:theorem:4}
\|e_i(t)\|_{\varepsilon}\leq e^{-\varphi_{\min} t}
\|e_i(0)\|_{\varepsilon} + \frac{D_\theta}{\varphi_{\min}}
\|\bar{\theta}_i-\hat{\theta}_i(\tau)\|_{\infty,[0,t]}, \
\bar{\theta}_i\in E_i(\theta_i)\cap [a,b].
\end{equation}

Let us now consider equations (\ref{eq:convergence_prototype:1}),
(\ref{eq:convergence_prototype:2}). We pick up a point $x'$, $y'$
which satisfies the following condition:
\begin{equation}\label{eq:theorem:4.5}
{x'}^2+{y'}^2=1.
\end{equation}
Solutions of (\ref{eq:convergence_prototype:2}) passing through this
point can be defined as follows:
\begin{equation}\label{eq:theorem:5}
\begin{split}
x_i(t,x',y')&=\cos\left(\int_0^t \gamma
\|\hat{s}_i(\tau)-s(\tau)\|_{\varepsilon} d\tau + \nu_x\right), \
x'=\cos(\nu_x), \ \nu_x\in[0,2\pi]\\
y_i(t,x',y')&=\sin\left(\int_0^t \gamma
\|\hat{s}_i(\tau)-s(\tau)\|_{\varepsilon} d\tau + \nu_y\right), \
y'=\sin(\nu_y), \ \nu_y\in[0,2\pi]
\end{split}
\end{equation}
This can be easily verified when writing
(\ref{eq:convergence_prototype:2}) in the system of polar
coordinates: $x_i=r \cos (\nu)$, $y_i=r \sin(\nu)$
\shortcite{Guckenheimer:2002}:
\begin{equation}\label{eq:theorem:5.5}
\begin{split}
\dot{r}&=\gamma \|\hat{s}_i-s\|_{\varepsilon} \cdot r(1-r)\\
\dot{\nu}&=\gamma \|\hat{s}_i-s\|_{\varepsilon}\cdot 1
\end{split}
\end{equation}
Given that $\bar{\theta}_i$ belongs to the interval $[a,b]$, there
is a number $\bar{h}(\bar{\theta}_i)\in[0,\pi]$ such that for all
$k\in\Numbers$ the following equivalence holds
\begin{equation}\label{eq:theorem:5.75}
\bar{\theta}_i= a+ \frac{b-a}{2}
\left(\cos(\bar{h}(\bar{\theta}_i)+2\pi k)+1\right).
\end{equation}
Hence according to (\ref{eq:convergence_prototype:1}),
(\ref{eq:theorem:5}) the norm
$\|\bar{\theta}_i-\hat{\theta}_i(\tau)\|_{\infty,[0,t]}$ can be
estimated from above as follows:
\begin{equation}\label{eq:theorem:6}
\|\bar{\theta}_i-\hat{\theta}_i(\tau)\|_{\infty,[0,t]}\leq\frac{b-a}{2}
\|\bar{h}(\bar{\theta}_i)-\nu_x+2 \pi k - \int_0^t \gamma
\|\hat{s}_i(\tau)-s(\tau)\|_{\varepsilon} d\tau\|_{\infty,[0,t]}
\end{equation}
Denoting
\[
\begin{split}
c&=\frac{D_\theta}{\varphi_{\min}}\frac{b-a}{2}; \ \
h(t,\bar{\theta}_i,k)=\bar{h}(\bar{\theta}_i)-\nu_x+2 \pi k -
\int_0^t \gamma \|\hat{s}_i(\tau)-s(\tau)\|_{\varepsilon} d\tau
\end{split}
\]
and taking into account (\ref{eq:theorem:4}), (\ref{eq:theorem:6})
we can conclude that the following holds along the solutions of
(\ref{eq:convergence_prototype})--(\ref{eq:convergence_prototype:2}):
\begin{equation}\label{eq:theorem:7}
\begin{split}
\|e_i(t)\|_{\varepsilon}&\leq e^{-\varphi_{\min} t}
\|e_i(0)\|_{\varepsilon} + c \|h(\tau,\bar{\theta}_i,k)\|_{\infty,[0,t]};\\
h(0,\bar{\theta}_i,k)&-h(t,\bar{\theta}_i,k)  = \int_0^t \gamma
\|e_i(\tau)\|_{\varepsilon} d\tau
\end{split}
\end{equation}
According to \shortcite{ArXive:Non-uniform:2006} (Theorem 1 and
Corollaries 2, 3) there exist $\gamma^\ast\in\Real_{>0}$ and
$h^\ast$ such that for a given bounded $e_i(0)$, all
$\gamma\in\Real_{>0}$, $\gamma<\gamma^\ast$ and
$h(0,\bar{\theta}_i,k)\geq h^\ast$ the norm
$\|e_i(\tau)\|_{\infty,[0,\infty]}$ is bounded and
\begin{equation}\label{eq:theorem:7.5}
\lim_{t\rightarrow\infty}h(t,\bar{\theta}_i,k)\in
[0,h(0,\bar{\theta}_i,k)].
\end{equation}
The value of $\gamma^\ast$, according to Corollary 3 in
\cite{ArXive:Non-uniform:2006}, can be determined from the following
inequality
\begin{equation}\label{eq:convergence_condition:1}
0<\gamma^\ast< \frac{\varphi_{\min}}{c}
\left(\ln\left(\frac{\kappa}{d}\right)\frac{\kappa}{\kappa-1}\left(2+\frac{\kappa}{1-d}\right)\right)^{-1},
\ \kappa\in \Real_{>1}, \ d\in(0,1)\subset\Real.
\end{equation}
The value of $h^\ast$ can be estimated from:
\begin{equation}\label{eq:convergence_condition:2}
\|e_i(t_0)\|_{\varepsilon}\leq
\left(\frac{\varphi_{\min}}{\gamma^\ast}\left(\ln
\frac{\kappa}{d}\right)^{-1}\frac{\kappa
-1}{\kappa}-c\left(2+\frac{\kappa}{1-d}\right)\right) h^\ast
\end{equation}
Given that $\|e_i(t_0)\|_{\varepsilon}$ in
(\ref{eq:convergence_condition:2}) is bounded from above for all
$t_0\geq 0$, $\|e_i(t_0)\|_{\varepsilon}\leq s_{\max} - s_{\min} +
D_\theta/\varphi_{\min}(b-a) $, condition
\begin{equation}\label{eq:convergence_condition:3}
h^\ast \geq \left((s_{\max}-s_{\min})+\frac{D_\theta
(b-a)}{\varphi_{\min}}\right)
\left(\frac{\varphi_{\min}}{\gamma^\ast}\left(\ln
\frac{\kappa}{d}\right)^{-1}\frac{\kappa
-1}{\kappa}-c\left(2+\frac{\kappa}{1-d}\right)\right)^{-1}
\end{equation}
together with (\ref{eq:convergence_condition:1}) imply that for all
$\hat{s}_i(t_0)\in\Omega_s$ and $h(0,\bar{\theta}_i,k)\geq h^\ast$
the norm $\|e_i(\tau)\|_{\infty,[0,\infty]}$ is bounded and property
(\ref{eq:theorem:7.5}) holds.

Notice that in the definition of $h(0,\bar{\theta}_i,k)$:
\begin{equation}\label{eq:theorem:7.25}
h(0,\bar{\theta}_i,k)=\bar{h}(\bar{\theta}_i)-\nu_x+2 \pi k
\end{equation}
the value of $k$ can be chosen arbitrarily large. Moreover,
$\bar{h}(\bar{\theta}_i)\in[0,\pi]$ for all
$\bar{\theta}_i\in[a,b]$. This implies that there exists a finite
$k'$ such that condition $h(0,\bar{\theta}_i,k')\geq h^\ast$ will be
satisfied for any fixed $h^\ast$ (i.e. for all $\gamma^\ast$
satisfying (\ref{eq:convergence_condition:1})) and all
$\bar{\theta}_i\in[a,b]$. In addition, the following will hold:
\begin{equation}\label{eq:theorem:7.75}
\lim_{t\rightarrow\infty}h(t,\bar{\theta}_i,k')\in
[0,h(0,\bar{\theta}_i,k')]\subset [0, \pi - \nu_x + 2\pi k'] \ \
\forall \in \bar{\theta}_i\in[a,b].
\end{equation}
Taking (\ref{eq:theorem:5}) into account we can conclude that
solutions $x_i(t,x',y')$ converge to a point in the interval
$[-1,1]$ as $t\rightarrow\infty$, and vector
$(x_i(t,x',y'),y_i(t,x',y'))$ makes no more than $k'$ full rotations
around the origin for all
$\theta_i\in[\theta_{\min},\theta_{\max}]$. Hence for a given
initial condition $x_i(0)=x'$, $y_i(0)=y'$, $\hat{s}_{i,0}\in
\Omega_s$  and ${\theta}_i\in[\theta_{\min},\theta_{\max}]$ the
estimate $\hat{\theta}_i(t)=a+(b-a)/2 \cdot (x_i(t,x',y')+1)$
converges to a point in $[a,b]$ as $t\rightarrow\infty$. We denote
this point by symbol $\hat{\theta}_i^\ast$.

Given that $\hat{\theta}_i(t)$ converges to a limit, there exists a
time instant $t^\ast$ such that for all $t\geq t^\ast$ the following
condition holds:
$|\hat{\theta}_i(t)-\hat{\theta}_i^\ast|<\mu_\infty$, where
$\mu_\infty \in\Real_{>0}$ is an arbitrarily small constant.
Therefore, taking condition (\ref{eq:Lipschitz}) into account, we
can conclude that for all $t\geq t^\ast$ derivative $\dot{e}_i$
satisfies the following equation:
\begin{equation}\label{eq:theorem:8}
\dot{e}_i=-\alpha(t) e_i + f_i(\xi(t),\theta_i) -
f_i(\xi(t),\hat{\theta}_i^\ast)+\mu_i(t) + \eta_i(t)
\end{equation}
where $|\mu_i(t)|\leq D_\theta \ \mu_\infty$ is a continuous
function.

Now we will show that the norm
$\normstwo{\theta_i}{E_i(\hat{\theta}_i^\ast)}$ can be bounded from
above by a $\mathcal{K}_\infty$-function  of $\Delta_\eta$. Consider
the term $f_i(\xi(t),\theta_i) - f_i(\xi(t),\hat{\theta}_i^\ast)$.
According to (\ref{eq:nonlinear_persistency}) there exists a
sequence of monotonically increasing time instances $t_j$,
$j=1,2,\dots$ such that $t_{j+1}-t_j\leq 2 T$ and
$|f_i(\xi(t_j),\theta_i) - f_i(\xi(t_j),\hat{\theta}_i^\ast)|\geq
\rho (\normstwo{\theta_i}{E_i({\hat{\theta}_i^\ast})})$.
Furthermore, according to (\ref{eq:Lipschitz_xi}),
(\ref{eq:xi_derivative}), the time-derivative of
$f_i(\xi(t),\theta_i) - f_i(\xi(t),\hat{\theta}_i^\ast)$ is bounded:
\[
\left|\frac{d}{dt} f_i(\xi(t),\theta_i) -
f_i(\xi(t),\hat{\theta}_i^\ast)\right|\leq 2 D_\xi \cdot \pd
\xi_\infty=D_f
\]
Hence the following estimate holds:
\begin{equation}\label{eq:theorem:9}
\begin{split}
 &\int_{t}^{t+L} |f_i(\xi(\tau),\theta_i) -
f_i(\xi(\tau),\hat{\theta}_i^\ast)| \geq \frac{\rho
(\normstwo{\theta_i}{E_i({\hat{\theta}_i^\ast})})^2}{2 D_f}\\
 & L = \max\left\{2 T, \frac{\rho(b-a)}{D_f}\right\}
\end{split}
\end{equation}

\noindent In order to proceed further we will need the following
lemma.

\begin{lem}\label{lem:filtered_pe} Consider the following
differential equation
\begin{equation}
\label{eq:filtered_pe_system} \dot{z}=-\varphi(t,z) + u(t) +
\eta(t), \ z_0=z(0)\in[z_{\min},z_{\max}]\subset \Real
\end{equation}
Let us suppose that

1) $\varphi(z)z\geq 0, \ \varphi_{\min}\leq {\pd \varphi(t,z)}/{\pd
z}\leq \varphi_{\max}$;

2) $u(t)\in L_{\infty}[0,\infty]\cap \mathcal{C}^1, \
\|u(t)\|_{\infty,[0,\infty]}\leq u_\infty$,
$\|\dot{u}(t)\|_{\infty,[0,\infty]}\leq \pd u_{\infty}$

3) $\eta(t)\in L_\infty[0,\infty]$,
$\|\eta(t)\|_{\infty,[0,\infty]}\leq \Delta$

4) there exist constants $L$, $\delta$ such that for all $t\geq 0$
\begin{equation}\label{eq:pe}
\int_{t}^{t+L}|u(\tau)|d\tau\geq \delta
\end{equation}

5) finally assume that the following inequality holds
\begin{equation}\label{eq:filtered_pe_condition}
\left(\frac{\delta}{L}\right)^2 - \Delta u_\infty > 0.
\end{equation}

Then for any $p\in\Real_{>0}$ there exist constants $L^\ast>0$ and
$\delta^\ast \geq (\left({\delta}/{L}\right)^2 - \Delta
u_\infty)/p$,
 such that
\begin{equation}\label{eq:filtered_pe_result}
\int_{t}^{t+L^\ast}|z(\tau)|d\tau\geq \delta^\ast\geq
\frac{1}{p}\left(\frac{\delta^2}{L} - \Delta u_\infty L \right) \
\forall \ t\geq 0
\end{equation}
\end{lem}
{\it Proof of Lemma \ref{lem:filtered_pe}.} We prove the lemma along
the lines of an argument provided in \shortcite{ArXive:Loria:2003}
(Property 1). Consider the time-derivative of $z u$:
\begin{equation}\label{eq:lemma_pe:1}
\frac{d}{dt}\left(z u\right) = (-\varphi(t,z) + u +\eta) u + z
\dot{u} \geq u^2 - |z|\left(\varphi_{\max}+\pd u_\infty\right) -
|u|\Delta
\end{equation}
According to (\ref{eq:lemma_pe:1}) for all  $t,t_0\in\Real_{\geq
0}$, $t\geq t_0$ the following inequality holds:
\begin{equation}\label{eq:lemma_pe:2}
z(t)u(t)-z(t_0)u(t_0)\geq \int_{t_0}^t u^2(\tau)d\tau -
\left(\varphi_{\max}+\pd u_\infty\right)\int_{t_0}^t |z(\tau)|d\tau
- \Delta \int_{t_0}^t|u(\tau)|d\tau
\end{equation}
Rearranging terms in (\ref{eq:lemma_pe:2}) yields
\[
\begin{split}
& \left(\varphi_{\max}+\pd u_\infty\right)\int_{t_0}^t
|z(\tau)|d\tau\geq z(t_0)u(t_0)-z(t)u(t) + \int_{t_0}^t
u^2(\tau)d\tau -  \Delta \int_{t_0}^t|u(\tau)|d\tau
\end{split}
\]
Notice that $z(t_0)u(t_0)-z(t)u(t)$ is bounded from below for all
$t\geq 0$. We denote this bound by symbol $M$. Furthermore,
according to the Holder inequality and property (\ref{eq:pe}), the
following estimate holds for all $t\geq 0$:
\[
\frac{\delta^2}{L}\leq\frac{1}{L} \left(\int_{t}^{t+L}|u(\tau)|d\tau
\right)^2 \leq \int_{t}^{t+L} u^2(\tau)d\tau.
\]
Hence for all time instances $t$: $(n+1) L \geq t-t_0\geq n L$,
where $n$ is a positive integer, we have
\begin{equation}\label{eq:lemma_pe:3}
\begin{split}
& \left(\varphi_{\max}+\pd u_\infty\right)\int_{t_0}^t
|z(\tau)|d\tau\geq M + n \frac{\delta^2}{L} -  \Delta
\int_{t_0}^t|u(\tau)|d\tau\\
& \geq M + n \frac{\delta^2}{L} - (n+1)\Delta u_\infty =  (M -\Delta
u_\infty L) + n \left(\frac{\delta^2}{L} - \Delta u_\infty L\right)
\end{split}
\end{equation}
According to the  requirements of the lemma, inequality
(\ref{eq:filtered_pe_condition}), the difference ${\delta^2}/{L} -
\Delta u_\infty L>0$ is a positive constant. Therefore, there
 exists $n=n'$ such that the right-hand side of
(\ref{eq:lemma_pe:3}) exceeds some $\delta'= ({\delta^2}/{L} -
\Delta u_\infty L)/p' \in\Real_{>0}$, $p'\in\Real_{>0}$. Choosing
$t'=\min_{t} \{t-t_0\}\geq n' L$ we can conclude that
\begin{equation}\label{eq:lemma_pe:4}
\left(\varphi_{\max}+\pd u_\infty\right)\int_{t_0}^{t'}
|z(\tau)|d\tau \geq \delta'
\end{equation}
Given that we could chose the value of $t_0$ arbitrarily in the
domain $\Real_{\geq 0}$, inequality (\ref{eq:lemma_pe:4}) is
equivalent to
\[
\int_{t}^{t+L^\ast} |z(\tau)|d\tau \geq \delta^\ast,
\]
where $L^\ast = t'-t_0$,
$\delta^\ast=\delta'/\left(\varphi_{\max}+\pd
u_\infty\right)=({\delta^2}/{L} - \Delta u_\infty L)/p$,
$p=p'\left(\varphi_{\max}+\pd u_\infty\right)$. {\it The lemma is
proven.}

Denoting $f_i(\xi(t),\theta_i) -
f_i(\xi(t),\hat{\theta}_i^\ast)=u(t)$, $\eta_i(t)+\mu_i(t)=\eta(t)$
we can observe that equation (\ref{eq:theorem:8}) is of the same
class as (\ref{eq:filtered_pe_system}) in the formulation of Lemma
\ref{lem:filtered_pe}. Furthermore, the following inequalities hold:
\begin{equation}\label{eq:theorem:10}
\Delta\leq \Delta_{\eta} + D_\theta \ \mu_\infty; \  \
\|u(t)\|_{\infty,[0,\infty]} \leq D_\theta
\normstwo{\theta_i}{E_i(\hat{\theta}_i^\ast)}\leq D_\theta (b-a)
\end{equation}
Notice that the value of $\mu_\infty$ in (\ref{eq:theorem:10}) can
be made arbitrarily small because $\hat{\theta}_i(t)$ converges to a
limit, and $\hat{\theta}_i^\ast$ can be chosen from its arbitrarily
small vicinity. Let us therefore chose $\hat{\theta}_i^\ast$ such
that $D_\theta \mu_\infty\leq \Delta_\eta$. Hence, in accordance
with Lemma \ref{lem:filtered_pe},  condition
\begin{equation}\label{eq:theorem:11}
\left(\frac{\rho^2(\normstwo{\theta_i}{E_i(\hat{\theta}_i^\ast)})}{2
D_f L}\right)^2 >  2 \Delta_\eta D_\theta (b-a)
\end{equation}
implies existence of constants $L^\ast$, $p\in\Real_{>0}$ such that
\begin{equation}\label{eq:theorem:12}
\begin{split}
& \int_{t}^{t+L^\ast} |e_i(\tau)|d\tau \geq \frac{1}{p}
\left(\left(\frac{\rho^2(\normstwo{\theta_i}{E_i(\hat{\theta}_i^\ast)})}{2
D_f}\right)^2\frac{1}{L}- \Delta u_\infty L \right)=\delta^\ast>0 \
\ \forall t\geq t^\ast.
\end{split}
\end{equation}

We are going to show now that the norm
$\normstwo{\theta_i}{E_i(\hat{\theta}_i^\ast)}$ is bounded from
above by a function
$\varepsilon_\theta(\Delta_\eta)\in\mathcal{K}_\infty$ for all
sufficiently small $\Delta_\eta$. Let us parameterize $\Delta_\eta$
as follows:
\begin{equation}\label{eq:theorem:12.5}
\Delta_\eta = \left(\frac{\rho^2(\varepsilon^\ast)}{2 D_f
L}\right)^2 \frac{1}{2 D_\theta (b-a)}, \
\varepsilon^\ast\in\Real_{> 0}.
\end{equation}
Parametrization (\ref{eq:theorem:12.5}) is always possible because
$\rho(\cdot)\in\mathcal{K}_\infty$. For all
$\normstwo{\theta_i}{E_i(\hat{\theta}_i^\ast)}> \varepsilon^\ast$
condition (\ref{eq:theorem:11}) is satisfied. Hence, according to
Lemma \ref{lem:filtered_pe} there exist constants $L^\ast$, $p$ such
that inequality (\ref{eq:theorem:12}) holds. Given that
$\delta^\ast, \ L^\ast, \ \varphi_{\min} \in\Real_{>0}$ there will
always exist a number $\Delta_\eta^\ast\in\Real_{>0}$ such that
$\Delta_\eta^\ast< (L^\ast)^{-1} \delta^\ast \varphi_{\min}/2$. This
implies that for all $\Delta_\eta\leq \Delta_\eta^\ast$ the
following inequality holds
\begin{equation}\label{eq:theorem:13}
\begin{split}
& \int_{t}^{t+L^\ast} \|e_i(\tau)\|_{\varepsilon}d\tau \geq
\frac{\delta^\ast}{2}, \
 \varepsilon=\frac{\Delta_\eta}{\varphi_{\min}}.
\end{split}
\end{equation}

Let us suppose that the norm
$\normstwo{\theta_i}{E_i(\hat{\theta}_i^\ast)}$ is greater than
$\varepsilon^\ast$.  In this case (\ref{eq:theorem:11}),
(\ref{eq:theorem:13}) hold and the integral
\begin{equation}\label{eq:theorem:14}
\int_{t^\ast}^t \|e_i(\tau)\|_{\varepsilon}d\tau
\end{equation}
grows unboundedly with $t$. On the other hand, according to
(\ref{eq:theorem:7}), (\ref{eq:theorem:7.5}) integral
(\ref{eq:theorem:14}) is bounded. Hence we have reached a
contradiction. This implies that
$\normstwo{\theta_i}{E_i(\hat{\theta}_i^\ast)}\leq\varepsilon^\ast$.
Given that $\rho(\cdot)\in \mathcal{K}_\infty$, the inverse
$\rho^{-1}(\cdot)$ is well defined and is a
$\mathcal{K}_\infty$-function. Therefore, taking
(\ref{eq:theorem:12.5}) into account, we can conclude that the
latter inequality is equivalent to:
\begin{equation}\label{eq:theorem:15}
\normstwo{\theta_i}{E_i(\hat{\theta}_i^\ast)} \leq \rho^{-1}\left(
 \left(8 \Delta_\eta D_\theta (b-a) D_f^2 L^2\right)^{1/4} \right)
\end{equation}
Thus we have just shown that there exists a point $x'$, $y'$ in
system
(\ref{eq:convergence_prototype})--(\ref{eq:convergence_prototype:2}),
and  parameters $\gamma$ and $\varepsilon$ such that the system
trajectories starting from this point converge into a small
neighborhood of $E_i(\theta_i)$ in finite time for all
$s_{i,0}\in\Omega_s$ and any given
$\theta_i\in[\theta_{\min},\theta_{\max}]$. The size of this
neighborhood can be characterized by a $\mathcal{K}_\infty$-function
of $\Delta_{\eta}$, when $\Delta_\eta$ is sufficiently small. Let us
now show that this convergence is uniform with respect to
$\theta_i$.

{\it 2. Uniformity.} Consider equation (\ref{eq:theorem:7.75}).
According to (\ref{eq:theorem:7}), (\ref{eq:theorem:7.75})
trajectories passing through a point $(x',y')$ satisfying
(\ref{eq:theorem:4.5}) at $t=0$ also satisfy the following
constraint:
\begin{equation}\label{eq:theorem:16}
\exists k'\in \Numbers: \ \
h(0)-h(\infty)=\gamma\int_{0}^\infty\|e_i(\tau,e_i(0),\theta_i,\eta_i(\tau))\|_\varepsilon
d\tau \leq {\pi - \nu_x + 2\pi k'}<\infty \
\end{equation}
for all $\theta_{i}\in[\theta_{\min},\theta_{\max}]$ and $e_i(0)$.
We will use this property to demonstrate that there is a point $(x',
y')$, $\sqrt{x'^2+y'^2} = 1$,
$\|\hat{\theta}_i(x')\|_{E_i(\theta_i)}\geq \Delta_0$,
$\Delta_0\in\Real_{>0}$,  such that   for any
$\theta_i\in[\theta_{\min},\theta_{\max}]$   the estimate
$\hat{\theta}_i(x_i(t,x',y'))$ converges into a set
\begin{equation}\label{eq:theorem:17}
\normstwo{\theta_i}{E_i(\hat{\theta}_i)}\leq \rho^{-1}\left(
 \left(8 \Delta_\eta D_\theta (b-a) D_f^2 L^2\right)^{1/4} \right)
\end{equation}
in finite time $T'(\theta_i)$ for all $t_0$,
$\hat{s}_{i,0}\in\Omega_s$, and stays there for all $t\geq t_0 +
T'(\theta_i)$. Furthermore, the value of $T'(\theta_i)$ is bounded
from above for all $\theta_{i}\in[\theta_{\min},\theta_{\max}]$. In
other words, there exists $T'_{\max}\in\Real_{>0}$:
\begin{equation}\label{eq:theorem:18}
T'(\theta_i)\leq T'_{\max} \ \forall \
\theta_i\in[\theta_{\min},\theta_{\max}].
\end{equation}
The fact that estimate $\hat{\theta}_i$ converges into a set
specified by (\ref{eq:theorem:17}) in finite time $T'(\theta_i)$ and
stays there for $t\geq t_0 + T'(\theta_i)$  for all $x', \ y': \
\sqrt{x'^2+y'^2}=1$ follows immediately from (\ref{eq:theorem:15}).
We must show, however, that (\ref{eq:theorem:18}) holds.

According to (\ref{eq:non-degeneracy_recognition}),
(\ref{eq:theorem:prototype_parameters}) there is a point
$\theta_0\in[a,b]$ such that $\|\theta_{0}\|_{E_i(\theta_i)}\geq
\Delta_\theta$ for every $\theta_i\in \Omega_\theta$. Hence, there
exists a point $\theta_{i,1}\in[a,b]$ such that
\[
\inf_{\bar{\theta}_i\in E_i(\theta_i)\cap[a,b]}\|
\bar{\theta}_i-\theta_{i,1}\|= \Delta_\theta
\]
Without loss of generality, suppose that the set
$\Omega_1=\{\bar{\theta}_i\in E_i(\theta_i)\cap [a,b]| \
\theta_{i,1}>\bar{\theta}_i \}$  is not empty\footnote{If $\Omega_1$
is empty then $\Omega_2=\{\bar{\theta}_i\in E_i(\theta_i)\cap [a,b]|
\ \theta_{i,1}<\bar{\theta}_i\}$ is not empty. We can proceed with
the same argument replacing interval $[0,\pi]$ with $[\pi, 2\pi]$
and $\sup$ with $\inf$ when appropriate.}. By symbol
$\theta_{i,\max}$ we denote $\theta_{i,\max}=\sup \{\Omega_1\}$. Let
us pick a point $\theta_{i,2}\in[a,b]$ according to the following
constraints
\begin{equation}\label{eq:theorem:18.03}
\begin{split}
|\theta_{i,2}-\theta_{i,1}|&=|\theta_{i,2}-\theta_{i,\max}|=\Delta_{\theta}/2,
\\
& \theta_{i,1}>\theta_{i,2}>\theta_{i,\max},
\end{split}
\end{equation}
and choose the value of $\nu_x$ in (\ref{eq:theorem:5}) such that
\[
\theta_{i,2}=a + \frac{b-a}{2}(\cos(\nu_x)+1), \ \nu_x\in[0,\pi].
\]
According to (\ref{eq:theorem:5.75}) there exist
$\bar{h}(\theta_{i,\max})$, $k$ such that
\[
\begin{split}
\theta_{i,\max}&=a +
\frac{b-a}{2}(\cos(\bar{h}(\theta_{i,\max})+2\pi k)+1), \
\bar{h}(\theta_{i,\max})\in[0,\pi], \ k\in\Natural.
\end{split}
\]
Given that $\theta_{i,2}>\theta_{i,\max}$ we set the value of $k=0$
and chose $\bar{h}(\theta_{i,\max})$ in accordance with the
following inequality:
\begin{equation}\label{eq:theorem:18.0625}
\nu_x < \bar{h}(\theta_{i,\max}).
\end{equation}
Because $
|\hat{\theta}_i(\cos(\nu_x))-\hat{\theta}_i(\cos(\nu_x'))|\leq
\frac{b-a}{2}|\nu_x-\nu_x'|$ for all $\nu_x,\nu_x'\in\Real$,
conditions (\ref{eq:theorem:18.03}), (\ref{eq:theorem:18.0625})
ensure existence of a constant $\nu_x'\leq
\bar{h}(\theta_{i,\max})$, $\nu_x'=\nu_x+\Delta_\theta/(2(b-a))$
such that
\begin{equation}\label{eq:theorem:18.07}
|\hat{\theta}_i(\cos(\nu_x))-\hat{\theta}_i(\cos(\nu_x''))|\leq
\Delta_{\theta}/4 \ \ \forall \ \nu_x''\in[\nu_x,\nu_x'].
\end{equation}
Hence,
\[
\|\hat{\theta}_i(\cos(\nu_x''))\|_{E_i(\theta_i)}\geq
\frac{\Delta_{\theta}}{4} \ \ \forall \ \nu_x''\in[\nu_x,\nu_x'].
\]
The inequality above implies that the values of
$\hat{\theta}_i(\cos(\nu_x''))$ are outside of the
$\Delta_\theta/4$-neighborhood of $E_i(\theta_i)$ for all
$\nu_x''\in[\nu_x,\nu_x']$. Furthermore, because
$\hat{\theta}_i(\cos(\cdot))$ is monotone (non-increasing) over
$[\nu_x,\bar{h}(\theta_{i,\max}))$, and
$\theta_{i,2}>\theta_{i,\max}$, there are no values of
$\nu_x''\in[\nu_x,\bar{h}(\theta_{i,\max}))$ such that
$\|\hat{\theta}_i(\cos(\nu_x''))\|_{E_i(\theta_i)}=0$.

Let us consider solutions of system
(\ref{eq:convergence_prototype})--(\ref{eq:convergence_prototype:2})
passing through the following point $x_i(0)=\cos(\nu_x)$,
$y_i(0)=\sin(\nu_x)$, $\hat{s}_i(0)\in\Omega_s$. Suppose that
$0<\gamma<\gamma^\ast$, and $\gamma^\ast$ satisfies
(\ref{eq:convergence_condition:3}) with
$h^\ast=\Delta_{\theta}/(2(b-a))$. Then, according to
\shortcite{ArXive:Non-uniform:2006} the sum $\nu_x + \gamma
\int_{0}^t \|e_i(\tau)\|_\varepsilon d\tau$ converges to a point in
$[\nu_x,\bar{h}(\theta_{i,\max})]$. Taking monotonicity and
continuity of function $\hat{\theta}_i(\cos(\nu_x''))$ for
$\nu_x''\in[\nu_x,\bar{h}(\theta_{i,\max})]$ into account, we can
conclude that trajectory $\hat{\theta}_i(x_i(t,x'(\theta_i)))$
enters the $\varepsilon^\ast$-neighborhood of $\theta_{i,\max}$ only
once for all $t\in[0,\infty]$.

Let us show that amount of time required for the system to enter
this neighborhood is bounded from above for all
$\theta_i\in\Omega_\theta$. Given that trajectory
$\hat{\theta}_i(x_i(t,x',y'))$ enters the
$\varepsilon^\ast$-neighborhood of $\theta_{i,\max}$ only once, we
shall show that the amount of time the system spends outside of
 this neighborhood is bounded from above for all
$\theta_i\in\Omega_\theta$. We prove this by contradiction. Suppose
that for any fixed $T'_0\in\Real_{>0}$ there is a
$\theta_i\in[\theta_{\min},\theta_{\max}]$ such that
$T'(\theta_i)\geq T'_0$. Consider dynamics of
(\ref{eq:convergence_prototype})--(\ref{eq:convergence_prototype:2})
when $s(t)=s_i(t,s_{i,0},\theta_i,\eta_{i}(t))$. Let us pick a
sequence of time instances $\{t_j\}_{j=1}^{\infty}$, such that
$t_{j+1}-t_{j}=D_T$, and $D_T\geq L^\ast$. For each  interval
$[t_j,t_{j+1}]$ we consider two possibilities:

1) the norm
$\|\hat{\theta}_i(t_j)-\hat{\theta}_i(\tau)\|_{\infty,[t_j,t_{j+1}]}\leq
\epsilon$, $\epsilon\in\Real_{>0}$, $\epsilon\leq D_\theta^{-1}
\Delta_\eta$, and

2) the norm
$\|\hat{\theta}_i(t_j)-\hat{\theta}_i(\tau)\|_{\infty,[t_j,t_{j+1}]}>
\epsilon$.

\noindent In case the first alternative applies, according to
(\ref{eq:theorem:13}) the following estimate holds
$\int_{t_j}^{t_{j+1}}\|e_i(\tau)\|_\varepsilon d\tau \geq
\delta^\ast$.  Hence $h(t_j)-h(t_{j+1})>\gamma \delta^\ast$. When
the second alternative holds, e.g.
$\|\hat{\theta}_i(t_i)-\hat{\theta}_i(\tau)\|_{\infty,[t_j,t_{j+1}]}
> \epsilon$, we can conclude, using inequality (\ref{eq:theorem:6}),
that
\[
\|\gamma
\int_{t_j}^{\tau}\|e_i(\tau_1)\|_{\varepsilon}d\tau_1\|_{\infty,[t_j,t_{j+1}]}
> \epsilon \frac{2}{b-a}.
\]
Given that $h(t)$ is monotone with respect to $t$ we obtain that
$h(t_j)-h(t_{j+1})>\epsilon 2/(b-a)$. Thus we have shown that
\[
h(t_j)-h(t_{j+1})>\min\{\gamma \delta^\ast, \epsilon
2/(b-a)\}=\Delta_h
\]
for all $j$ such that
$\normstwo{\hat{\theta}_i(\tau)}{E_i(\theta_i)}\geq
\varepsilon^\ast$ for all $\tau\in [t_j,t_{j+1}]$.  Given that
$h(t)$ is non-increasing and $T'$ is arbitrarily large, there will
be a time instance $t_m\leq T'$ such that $\sum_{j}^{m}
h(t_{j})-h(t_{j+1}) \geq m \Delta_h > {\pi - \nu_x + 2\pi k'}$.
This, however, contradicts to (\ref{eq:theorem:16}). Hence property
(\ref{eq:theorem:18}) is proven.

{\it 3. Structurally stable prototype.} So far we have shown that
for the given system
(\ref{eq:convergence_prototype})--(\ref{eq:convergence_prototype:2})
there exists a non-empty set of parameters $\gamma$, $\varepsilon$,
and  $x', y': \sqrt{{x'}^2+{y}'^2}=1$ such that trajectories
$x_i(t,x',y')$, $y_i(t,x',y')$ converge to a point on the unit
circle in $\Real^2$, and variable $\hat{\theta}_i(x_i(t,x',y'))$
reaches a given small vicinity of $E_i(\theta_i)$ (see
(\ref{eq:theorem:17})) within finite time $T'_{\max}$ for all
$\theta_i\in[\theta_{\min},\theta_{\max}]$.

Let us now consider perturbed system
(\ref{eq:convergence_prototype_struct_stable}),
(\ref{eq:convergence_prototype_struct_stable:2}) where
$\delta\in\Real_{>0}$ and initial conditions are selected in a
neighborhood of $x'$, $y'$:
\begin{equation}\label{eq:theorem:18.25}
(x_i(0),y_i(0))\in\Omega(x',y')=\{(x,y)\in\Real^2|
\sqrt{(x-x')^2+(y-y')^2}\leq \delta_r \}, \ \delta_r\in\Real_{>0}.
\end{equation}
In order to distinguish solutions of
(\ref{eq:convergence_prototype_struct_stable}),
(\ref{eq:convergence_prototype_struct_stable:2}) from  the solutions
of unperturbed system
(\ref{eq:convergence_prototype})--(\ref{eq:convergence_prototype:2}),
we denote the latter by symbols $x_i^\ast(t,x_i(0),y_i(0))$,
$y_i^\ast(t,x_i(0),y_i(0))$, and
$\hat{s}_i^\ast(t,\theta_i,s_{i,0},\eta_i(t))$. For the sake of
notational compactness we also denote the state vector of the $i$-th
subsystem of
(\ref{eq:convergence_prototype})--(\ref{eq:convergence_prototype:2})
as $\bfq_i^\ast=(\hat{s}_i^\ast,x_i^\ast,y_i^\ast)$, and the state
vector of the $i$-th subsystem of
(\ref{eq:convergence_prototype_struct_stable}),
(\ref{eq:convergence_prototype_struct_stable:2}) as $\bfq_i$.

Solutions of (\ref{eq:convergence_prototype_struct_stable}),
(\ref{eq:convergence_prototype_struct_stable:2}) are bounded:
\begin{equation}\label{eq:theorem:18.5}
\begin{split}
\|\hat{s}_i(t,\hat{s}_{i,0},\eta_i(t))\|_{\infty,[0,\infty]}&\leq
|\hat{s}_{i,0}|+ (\max\{|a|,|b|\}
D_{\theta}+\Delta_\eta)/\varphi_{\min},\\
\|x_i(t,x_i(0),y_i(0))\|_{\infty,[0,\infty]}&\leq
\max\{1,\sqrt{x_i(0)^2+y_i(0)^2}\},\\
\|y_i(t,x_i(0),y_i(0))\|_{\infty,[0,\infty]}&\leq
\max\{1,\sqrt{x_i(0)^2+y_i(0)^2}\}.
\end{split}
\end{equation}
Hence for all $\hat{s}_{i}(0), x_i(0), y_i(0)\in
\Omega_s\times\Omega(x',y')$ there exists a constant $D_0$ such that
$\|\bfq_i(t)\|_{\infty,[0,\infty]}\leq D_0$ for all $\theta_i$. Let
us rewrite (\ref{eq:convergence_prototype_struct_stable}),
(\ref{eq:convergence_prototype_struct_stable:2}) as follows:
\begin{equation}\label{eq:theorem:18.75}
\begin{split}
\dot{\hat{s}}_i&=-\varphi_i(\hat{s}_i)+f_i(\xi(t),\hat{\theta}_i(x_i))\\
\dot{x}_i&=\gamma \|\hat{s}_i-s\|_{\varepsilon} \left(x_i-y_i-x_i(x_i^2+y_i^2)\right) +
\gamma \delta \cdot \varepsilon_x(x_i,y_i) \\
\dot{y}_i&=\gamma
\|\hat{s}_i-s\|_{\varepsilon}\left(x_i+y_i-y_i(x_i^2+y_i^2)\right) +
\gamma \delta \cdot \varepsilon_y(x_i,y_i),
\end{split}
\end{equation}
where
\[
\begin{split}
\varepsilon_{x}(x_i(t),y_i(t))&=x_i(t)-y_i(t)-x_i(t)(x_i^2(t)+y_i^2(t)); \\
\varepsilon_{y}(x_i(t),y_i(t))&=x_i(t)+y_i(t)-y_i(t)(x_i^2(t)+y_i^2(t))
\end{split}
\]
The right-hand side of
(\ref{eq:convergence_prototype})--(\ref{eq:convergence_prototype:2})
is locally Lipschitz in $\hat{s}_i$, $x_i$, $y_i$ (and so is the
right-hand side of (\ref{eq:convergence_prototype_struct_stable}),
(\ref{eq:convergence_prototype_struct_stable:2})). We denote its
corresponding  Lipschitz constant in the domain specified by
(\ref{eq:theorem:18.5}) by symbol $L_i(D_0)$. Furthermore, provided
that (\ref{eq:theorem:18.5}) holds,  $\varepsilon_x(x_i(t),y_i(t))$,
$\varepsilon_y(x_i(t),y_i(t))$ are globally bounded with respect to
$t$. Let us denote this bound by symbol $B$:
\[
\max \left\{ \|\varepsilon_x(x_i(t),y_i(t))\|_{\infty,[0,\infty]},
\|\varepsilon_y(x_i(t),y_i(t))\|_{\infty,[0,\infty]} \right\} = B
\]
For the sake of notational compactness let us rewrite
(\ref{eq:theorem:18.75}) as follows:
\begin{equation}\label{eq:theorem:18.875}
\dot{\bfq}_i=\bff(\bfq_i,s(t),\xi(t))+ \gamma \delta \cdot
\bfg(\bfq_i),
\end{equation}
where $\bff(\bfq_i,s(t),\xi(t))$ and $\bfg(\bfq_i)$ are defined to
copy the right-hand side of (\ref{eq:theorem:18.75}). Notice that
$\|\bff(\bfq_i,s(t),\xi(t))\|\leq L_i(D_0) \|\bfq_i\|$,
$\|\bfg(\bfq_i)\|\leq B\sqrt{2}$.

According to the theorem on continuous dependence of solutions of an
ODE on parameters and initial conditions (see, for instance,
\shortcite{Khalil:2002} Theorem 3.4, page 96) the following holds:
\begin{equation}\label{eq:theorem:19}
\|\bfq_i(t)-\bfq_i^\ast(t)\|\leq \|\bfq_i(t_0)-\bfq_i^\ast
(t_0))\|e^{L_i(D_0)(t-t_0)} + \frac{\delta \gamma B
\sqrt{2}}{L_i(D_0)}\left(e^{L_i(D_0)(t-t_0)}-1\right).
\end{equation}
When the values of $\hat{s}_{i,0}$ and $\hat{s}_{i,0}^\ast$ coincide
estimate (\ref{eq:theorem:19}) implies that
\begin{equation}\label{eq:theorem:19.5}
\|\bfq_i(t)-\bfq_i^\ast(t)\|\leq \delta_r  e^{L_i(D_0)(t-t_0)} +
\frac{\delta \gamma
B\sqrt{2}}{L_i(D_0)}\left(e^{L_i(D_0)(t-t_0)}-1\right).
\end{equation}
This assures existence of $\delta_r\in\Real_{>0}$,
$\delta\in\Real_{>0}$ such that for a fixed, yet arbitrarily large,
time $T''(\delta_r,\delta)>T'_{\max}$ solutions of system
(\ref{eq:convergence_prototype_struct_stable}),
(\ref{eq:convergence_prototype_struct_stable:2}) passing through a
point from $\Omega(x',y')$ at $t=t_0$ will remain within a fixed,
yet arbitrarily small, neighborhood of a solution of system
(\ref{eq:convergence_prototype})--(\ref{eq:convergence_prototype:2})
with initial conditions $x_i(t_0)=x'$, $y_i(t_0)=y'$. The value of
$T_{\max}'$ does not depend on $\delta_r$, $\delta$.

Taking (\ref{eq:theorem:5.5}) into account, we can conclude that the
set $x_i^2+y_i^2=1$ is globally attracting in the state space of
system (\ref{eq:convergence_prototype_struct_stable}),
(\ref{eq:convergence_prototype_struct_stable:2}) for almost all
initial conditions (except when $x_i(t_0)=0$, $y_i(t_0)=0$). This
implies that solutions starting in $\Omega(x',y')$ will remain
there.  In addition, according to (\ref{eq:theorem:5}), for any
$t_0\geq 0$ a $\delta_r$-vicinity of $(x',y')$ will be visited
within at least time $t'\leq t_0 + 2\pi/(\gamma \cdot \delta)$.
Hence we have just shown that for all $t_0\geq0$ solutions starting
at $\Omega_s\times \Omega(x',y')$ approach the target set within a
fixed time $T_{\max}'$ and stay in its vicinity for arbitrarily long
time $T''(\delta_r,\delta)$. The latter time is a function of
$\delta_r$, $\delta$: the smaller the values of $\delta_r$,
$\delta$, the larger the value of $T''(\delta_r,\delta)$.

{\it 4. Realizability.} Let us finally show that system
(\ref{eq:convergence_prototype_struct_stable}),
(\ref{eq:convergence_prototype_struct_stable:2}) can be realized by
a recurrent neural network. More precisely, we wish to prove that
there exists a system (\ref{eq:RNN_recurrent}) such that
$\bfx=\zetavec_1\oplus\zetavec_2\oplus\cdots\oplus\zetavec_{N_f}$,
$\zetavec_i\in\Real^3$,
$\zetavec_{i}=\zeta_{i,1}\oplus\zeta_{i,2}\oplus\zeta_{i,3}$,
$i=\{1,\dots,N_f\}$ and solutions $\zetavec_i(t,\bfq_{i,0})$ are
sufficiently close to $\bfq_i(t,\bfq_{i,0})$, where $\bfq_{i,0}\in
\Omega_s\times\Omega(x',y')\subset\Real^3$.

It is clear that the right-hand side of
(\ref{eq:convergence_prototype_struct_stable}),
(\ref{eq:convergence_prototype_struct_stable:2}) is a continuous and
locally Lipschitz function. To proceed further we use the following
result by Cybenko \shortcite{Cybenko}:
\begin{thm}[Cybenko, 1989]\label{theorem:Cybenko} Let $\sigma:\Real\rightarrow\Real$ be
any continuous sigmoid-type function. Then finite sums of the form
\[
G(\zetavec)=\sum_{j=1}^N \alpha_j \sigma (\omegavec_j^T \zetavec +
\beta_j), \ \zetavec\in\Real^m, \ \omegavec_j\in\Real^m, \ \alpha_j,
\ \beta_j\in\Real
\]
are dense in $\mathcal{C}[0,1]^n$.
\end{thm}

According to Theorem \ref{theorem:Cybenko}, for any arbitrarily
small $\varepsilon_N\in\Real_{>0}$, any given bounded intervals
$\Omega_x\subset\Real$, $\Omega_y\subset\Real$, and any
\[
s(t),\xi(t): \ \max\{\|s(t)\|_{\infty,[0,\infty]},
\|\xi(t)\|_{\infty,[0,\infty]}\}< M, \ M\in\Real_{>0},
\]
there exist $N\in\Natural$, $\omegavec_{j}\in\Real^5$,
$\alpha_{j}\in\Real$, $\beta_{j}\in\Real$, $j=1,2,\dots,N$ such that
\begin{equation}\label{eq:theorem:20}
\left|\sum_{j=1}^N \alpha_j \sigma (\omegavec_j^T \cdot
(\xi(t)\oplus s(t)\oplus\zetavec_i ) + \beta_j) -
\bff(\zetavec_i,s(t),\xi(t)) - \gamma \delta \cdot
\bfg(\zetavec_i)\right|<\varepsilon_N,
\end{equation}
where $\zetavec_i\in\Omega_s\oplus\times\Omega_x\times\Omega_y$. It
follows from (\ref{eq:theorem:20}) that there exist $N$,
$\omegavec_j$, $\alpha_j$, $\beta_j$ such that
\begin{equation}\label{eq:theorem:20.5}
\sum_{j=1}^N \alpha_j \sigma (\omegavec_j^T \cdot (\xi(t)\oplus
s(t)\oplus\zetavec_i ) + \beta_j) = \bff(\zetavec_i,s(t),\xi(t)) +
\gamma \delta \cdot \bfg(\zetavec_i) +
\Delta(\zetavec_i,s(t),\xi(t)),
\end{equation}
where $\Delta(\zetavec_i,s(t),\xi(t))$ is continuous and
\[
|\Delta(\zetavec_i,s(t),\xi(t))| < \varepsilon_N.
\]

Let us chose $\Omega_x=[-v,v]$, $\Omega_y=[-v,v]$ where
$v\in\Real_{>0}$, $v > 1$ and consider the dynamics of
\begin{equation}\label{eq:theorem:21}
\dot{\zetavec}_i=\bff(\zetavec_i,s(t),\xi(t))+\gamma \delta \cdot
\bfg(\zetavec_i) + \Delta(\zetavec_i,s(t),\xi(t)).
\end{equation}
System (\ref{eq:theorem:21}) has a globally attracting invariant set
(for almost all initial conditions) which can be characterized as
follows
\[
\{\zetavec_i\in\Real^3| \ 1-\rho(\varepsilon_N)\leq
{\zeta_{i,2}^2+\zeta_{i,3}^2}\leq 1+\rho(\varepsilon_N)\}, \
\rho\in\mathcal{K}_\infty.
\]
This follows immediately from the fact that
(\ref{eq:theorem:18.875}) is structurally stable and has a globally
attracting invariant set (for almost all initial conditions).
Furthermore, for any given $\varepsilon_N$ and a bounded set of
initial conditions $\Omega_\zeta (r)=\{\zetavec_i\in\Real^3| \
\|\zetavec_i\|\leq r, \ r\in\Real_{>0}\}$ there exists constant
$B_1$ such that $\|\zetavec_i(t)\|_{\infty,[0,\infty]}<B_1$ . Hence
solutions of system
\begin{equation}\label{eq:theorem:22}
\dot{\zetavec}_i=\sum_{j=1}^N \alpha_j \sigma (\omegavec_j^T \cdot
(\zetavec_i\oplus s(t)\oplus \xi(t) + \beta_j)
\end{equation}
are bounded for all initial conditions from $\Omega_\zeta (r)$
provided that inequality (\ref{eq:theorem:20}) holds over
sufficiently large intervals $\Omega_x$, $\Omega_y$ (for
sufficiently large $v$). Furthermore, given that $\varepsilon_N$ is
sufficiently small, solutions of (\ref{eq:theorem:22}) enter domain
$\Omega_s\times\Omega(x',y')$ specified by (\ref{eq:theorem:18.25})
in finite time. Finally, according to equality
(\ref{eq:theorem:20.5}) and Theorem 3.4 in \shortcite{Khalil:2002},
solutions of (\ref{eq:theorem:22}) starting in $\Omega(x',y')$
satisfy the following inequality:
\begin{equation}\label{eq:theorem:23}
\|\bfq_i(t,\bfq_{i,0})-\zetavec_i(t,\bfq_{i,0})\|\leq
\frac{\varepsilon_N}{L_i(D_0)}\left(e^{L_i(D_0)(t-t_0)}-1\right), \
\bfq_{i,0}\in\Omega_s\times\Omega(x',y').
\end{equation}
Hence, for any $t\geq 0$, solutions of (\ref{eq:theorem:22})
starting from  $\Omega_\zeta (r)$ approach the target set within a
fixed time (dependant on $\delta$) and stay in its vicinity
arbitrary long provided that $\delta$ and $\varepsilon_N$ are
sufficiently small. The possibility of the latter follows from
Theorem \ref{theorem:Cybenko}.


Taking (\ref{eq:theorem:23}), (\ref{eq:theorem:19.5}),
(\ref{eq:theorem:1}), (\ref{eq:convergence_prototype_struct_stable})
into account we conclude the proof by choosing $h_{f,i}(\bfx,s)$,
$h_{\theta,i}(\bfx,s)$ as follows
\begin{equation}\label{eq:theorem:24}
\begin{split}
h_{f,i}(\bfx,s)&=h_{f,i}(\zetavec_1\oplus\cdots\oplus\zetavec_{N_f},s)=s-\zeta_{i,1},\\
h_{\theta,i}(\bfx,s)&=h_{\theta,i}(\zetavec_1\oplus\cdots\oplus\zetavec_{N_f},s)=a
+ \frac{b-a}{2}(\zeta_{i,2}+1)
\end{split}
\end{equation}
{\it The theorem is proven.}

\noindent Before concluding this section we would like to provide
several remarks regarding Theorem \ref{theorem:RNN_existence}.

\begin{rem}[Read-out from the outputs]\normalfont
As follows from the theorem the class of signal
$s(t)=s_i(t,s_{i,0},\theta_i,\eta_i(t))$, e.g. parameter $i$, can be
inferred from the values of  $h_{f,j}(\bfx(t),s(t))$,
$j=\{1,\dots,N_f\}$ within a finite interval of time. The values of
$h_{f,i}(\bfx(t),s(t))$ should approach a small neighborhood of zero
and stay there for a sufficiently long time. The estimate  of
$\theta_i$ up to its equivalence class is available from the values
of $h_{\theta,i}(\bfx(t),s(t))$ over the same interval.

From a practical viewpoint, however, it is preferable to read-out
from the RNN outputs explicitly, rather than having to satisfy
ourselves with the existence of two sets of read-out functions, for
state and input, respectively, of the RNN. Even though this option
is not stated explicitly in Theorem \ref{theorem:RNN_existence}, it
can be easily shown that the preferred option can, indeed, be
realized. Adding to recurrent subsystem (\ref{eq:RNN_recurrent}) a
{\it feed-forward} part realizing continuous "output" functions
(\ref{eq:theorem:24}) enables explicit read-out from the RNN
outputs.
\end{rem}

\begin{rem}[Convergence to an attractor]\normalfont Theorem
\ref{theorem:RNN_existence} does not imply that recognition of a
class of the input signal $s(t)$ involves convergence of the RNN
state to an attractor. Yet its formulation does not exclude this
option either. In fact, when $f_i(\xi(t),\theta_i)$ satisfies some
additional restrictions (e.g. linear or monotone parametrization
with respect to $\theta_i$), it is possible to replace
(\ref{eq:convergence_prototype:1}),
(\ref{eq:convergence_prototype:2}) with another prototype system:
one that converges to a point attractor exponentially
\shortcite{IEEE_TAC_2007}. This implies that it depends
substantially on the properties of $f_i(\xi(t),\theta_i)$ whether
the state of a network will behave intermittently or asymptotically
converge to an attractor. It is important, however, that in both
cases the recognition problem will be successfully solved by a RNN.
\end{rem}

\begin{rem}[Multidimensional uncertainty]\normalfont Even though the
theorem applies to the case where  $\theta_i$ is a scalar, it can be
trivially extended to the case where  uncertain parameters are
vectors from a bounded domain $\Omega_{\theta,d}\subset\Real^d$. To
do so one needs to find a Lipschitz mapping
$\lambdavec:\Real\rightarrow\Real^d$ such that for a given small
$\varepsilon_{\lambda}\in\Real_{>0}$ the following property holds:
\[
\forall \ \thetavec_i\in\Omega_{\theta,d} \ \exists  \
\theta_i\in\Omega_{\theta}: \
\|\thetavec_i-\lambdavec(\theta_i)\|<\varepsilon_\lambda
\]
Hence the problem will reduce to the scalar case to which Theorem
\ref{theorem:RNN_existence} applies.
\end{rem}


%
%
%

\section{Conclusion}\label{Conclusion}


We provided a theoretical justification to the important question
why an RNN with fixed weights can serve as a universal adaptive
classifier of both static and dynamic inputs.  In addition to
providing an existence proof we have proven that the number of
dynamical states in an RNN recognizing $n$ different signals
$s_i(t)$ can be as small as $3 n$, i.e. grows linearly with the size
of the set of uncertain signals to be classified.

We stated the classification and recognition problems in a
behavioral context in which, over time, the desired input-output
relationship is achieved. Finding a solution corresponds to a
network dynamics in which the state reaches a given neighborhood of
the a-priori specified set and stays there for sufficiently long
time, provided that input to the network belongs to a given class
(Problem \ref{problem:RNN_existence}).  With these ramifications,
RNN solve the problem of adaptively classifying time-dependent
signals. We did not set out to guarantee, however, that the state of
the RNN will asymptotically converge to an equilibrium or its small
vicinity as a result of recognition. On the other hand the amount of
time a network would spend in the vicinity of a target set can be
made sufficiently large to qualify as a practical solution to the
classification problem. For classification, after all, asymptotic
convergence is not needed.

In physics and nonlinear dynamics the phenomenon that the  state of
a system reaches a neighborhood of a set and stays there
sufficiently long, yet inevitably escaping -- only to get caught
again, is called (chaotic) itinerancy \shortcite{Chaos:Kaneko:2003};
the set is referred to as an attractor-ruin. These descriptive
concepts are currently recognized as a possible mathematical basis
for modeling brain activity \shortcite{Tsuda:1991,Tsuda:2004}. We
envisage that our current result supports this idea, by showing the
considerable power of these systems to perform adaptive
classification.

\bibliographystyle{apacite}
\bibliography{existence_of_adaptive_RNN_new}

\end{document}